\documentclass[11pt]{article}
\usepackage[margin=1in]{geometry}
\usepackage{amsmath,amssymb,amsthm,mathtools}
\usepackage{booktabs}
\usepackage{mathrsfs}
\usepackage[hidelinks]{hyperref}
\hypersetup{pdftitle={Exact Finite-Horizon Difference Prophet Inequalities for Randomized Single-Threshold Rules},pdfauthor={Igor Kleiner}}
\usepackage{microtype}

\newtheorem{theorem}{Theorem}[section]
\newtheorem{lemma}[theorem]{Lemma}
\newtheorem{corollary}[theorem]{Corollary}
\newtheorem{proposition}[theorem]{Proposition}
\theoremstyle{remark}
\newtheorem{remark}[theorem]{Remark}

\newcommand{\Ts}{\mathscr T_{\mathrm{s,r}}}
\newcommand{\Funit}{\mathcal F_{[0,1]}}

\title{Exact Finite-Horizon Difference Prophet Inequalities\\for Randomized Single-Threshold Rules}
\author{Igor Kleiner\\\small Holon Institute of Technology, Israel\\\small\texttt{igorkl@hit.ac.il}}
\date{}

\begin{document}
\maketitle

\begin{abstract}
Consider the sequential selection of one of several independent and
identically distributed observations from a known law on $[0,1]$. Before the final observation, a
single-threshold rule accepts the first value above the threshold and
randomizes at equality. If no such acceptance occurs, it takes the last
observation. Goldenshluger, Malinovsky and Zeevi express the sharp
worst-case difference between the expected maximum and the best expected
reward of such a rule as the value of an infinite zero-sum game. We solve
this game for every horizon of at least two observations. One scalar
equation determines a minimax-optimal quantile, and the unique
least-favorable distribution has three atoms. The sharp constants increase
strictly to $0.1736372853\ldots$. We derive their leading logarithmic
finite-horizon correction and give an explicit quantile rule that
guarantees the limiting bound at every horizon.
\end{abstract}



\section{Introduction}\label{sec:introduction}

One of the simplest ways to select an observation sequentially is to
accept the first value above a prescribed threshold. In a finite sequence,
however, the consequences of rejecting an observation depend on how many
observations remain. An optimal stopping rule may therefore use a different
threshold at each stage. A single-threshold rule dispenses with this
flexibility: once the threshold has been chosen, the same comparison is
used at each decision before the final observation. Such a rule is simple
to implement, but its expected reward can be smaller than that of a
time-dependent rule. We seek to determine how large the expected advantage
of a prophet who sees the entire sequence in advance can be, even when
the threshold is chosen optimally for the distribution of the observations.

We consider independent observations with a common known distribution
supported on $[0,1]$. Exactly one observation must be selected. The
threshold is fixed before the observations begin; ties at the threshold
are resolved by independent coins with a fixed acceptance probability.
If no earlier observation is accepted, the last one is taken. For a given
distribution, we compare the greatest expected reward obtainable by such a
rule with the expected maximum of the observations. Our object is the
largest difference between these two expectations over all distributions
on $[0,1]$, for each finite number of observations. Thus the distribution
may be used in choosing the threshold. The comparison is additive
rather than multiplicative: on the unit interval it measures a loss in
reward on a common scale. This is a comparison with full foresight, not
merely with the best unrestricted stopping rule.

Additive prophet inequalities for unrestricted stopping rules were studied
by Hill and Kertz~\cite{HK81,HK82}. Threshold comparisons and variants
with observation costs were considered by Samuel-Cahn~\cite{SC84,SC92}.
For the class of randomized single-threshold rules considered here,
Goldenshluger, Malinovsky and Zeevi~\cite{GMZ} characterize the sharp
additive bound as the value of an infinite zero-sum game. Their
Section~5.2 computes bounds on the game value by discretizing the strategy
interval and solving a linear program. Here we ask whether the difference
game itself can be solved without discretization. In particular, can one
identify a quantile that works for every underlying law, and exhibit a
law against which no single-threshold rule can do better?

We give an exact solution for every horizon of at least two observations.
A single scalar equation with a unique root determines the optimal
quantile and the sharp bound. The quantile depends only on the horizon;
the threshold itself is obtained from the known distribution. We also
construct a three-point distribution for which no admissible rule has a
smaller loss than the sharp bound, and prove that this least-favorable
law is unique. The sharp bounds increase strictly with the horizon and
converge to $0.1736372853\ldots$. We derive the leading logarithmic
correction to this limit and show that a simpler quantile rule guarantees
the limiting bound at every finite horizon.

The proof starts from the game representation of Goldenshluger, Malinovsky
and Zeevi. A change of variable makes the payoff convex in the decision
variable, so that one quantile suffices in place of a mixture. This
reduction is the single-choice specialization of the static-threshold
construction of Jiang, Ma and Zhang~\cite{JMZ}; we give a direct proof
for the difference kernel in Section~\ref{sec:convex}. The two branches of the kernel can
then be maximized separately. Their equality determines the optimal
quantile, and a suitable mixture of their maximizing points gives the
matching lower bound.

Section~\ref{sec:setting} specifies the rule class and recalls the payoff
kernel. Sections~\ref{sec:convex}--\ref{sec:saddle} solve the game.
Section~\ref{sec:stopping} describes the stopping rule and proves the
uniqueness of the least-favorable law. Sections~\ref{sec:asympt}
and~\ref{sec:uniform} give the asymptotic analysis, strict increase with
the horizon and the explicit uniform bound.

\section{Stopping rules and the difference kernel}\label{sec:setting}

The common distribution need not be continuous. At an atom, a threshold
alone does not specify the probability of acceptance, and it is this
probability that determines the chance of reaching a later observation.
We therefore include tie-breaking in the definition of the stopping rule
before passing to the quantile parametrization.

Let $n\ge2$ be the number of observations, and let $\Funit$ denote the class
of probability distributions supported on $[0,1]$. Let $X_1,\ldots,X_n$
be independent with a common known distribution $F\in\Funit$. We also use
$F$ for its cumulative distribution function:
$F(t)=\mathbb P_F(X_1\le t)$. A randomized single-threshold rule fixes
$\theta\in[0,1]$ and $p\in[0,1]$ before observation begins. At each time
$i<n$, it accepts $X_i>\theta$, rejects $X_i<\theta$, and accepts
$X_i=\theta$ with probability $p$. The acceptance coins are independent of
one another and of the observations. If all first $n-1$ observations are
rejected, the rule takes $X_n$ regardless of its value. Write
$\tau_{\theta,p}$ for the selected index and $\Ts$ for this class of rules.
This is the randomized threshold class of Goldenshluger, Malinovsky and
Zeevi (GMZ)~\cite{GMZ}.

The prophet sees all observations and selects their maximum. For a rule
$\tau\in\Ts$, define
\[
 M_n(F):=\mathbb E_F\max_{1\le i\le n}X_i,
 \qquad V_n(\tau;F):=\mathbb E_F X_\tau.
\]
The second expectation includes the acceptance coins. We call
$M_n(F)-V_n(\tau;F)$ the regret of the rule. The sharp worst-case regret is
\begin{equation}\label{eq:regretdef}
 v_n:=\sup_{F\in\Funit}\inf_{\tau\in\Ts}
       \{M_n(F)-V_n(\tau;F)\}.
\end{equation}
The order of the two operations in \eqref{eq:regretdef} reflects that $F$ is
known before the rule is chosen. Parametrizing a threshold by its quantile
allows rules for different laws to be compared in the same coordinate.

For $u\in[0,1]$, write $F(u-)=\lim_{t\uparrow u}F(t)$ and
$\Delta F(u)=F(u)-F(u-)$ for the mass at $u$, with $F(0-)=0$.
The probability of rejecting one observation under a threshold $\theta$
and tie probability $p$ is $F(\theta)-p\Delta F(\theta)$.
We use this rejection probability, denoted by $x$, to parametrize the rule.
It determines how likely the rule is to reach later observations and is
convenient even when $F$ has atoms or gaps in its support.

For $0<x<1$, define the left-continuous quantile
\[
 \theta_x=F^{\leftarrow}(x):=\inf\{u\in[0,1]:F(u)\ge x\}.
\]
When $\Delta F(\theta_x)>0$, set
\begin{equation}\label{eq:tie}
 p_x=\frac{F(\theta_x)-x}{\Delta F(\theta_x)}.
\end{equation}
The inequalities $F(\theta_x-)\le x\le F(\theta_x)$ show that
$p_x\in[0,1]$. When $\Delta F(\theta_x)=0$, one has $F(\theta_x)=x$;
ties then have probability zero, and we set $p_x=0$. In either case the probability of rejecting one
observation is
\[
 F(\theta_x)-p_x\Delta F(\theta_x)=x.
\]
Write $\tau_x^F$ for this rule. At $x=0$ it accepts the first observation;
at $x=1$ it waits until the last. The coordinate $x$ is both a quantile
level and the one-step rejection probability. Independence implies that
the rule reaches the final observation with probability $x^{n-1}$.
The following equivalence
also covers thresholds at the edge of a gap in the support of $F$.

\begin{lemma}[Quantile parametrization]\label{lem:parametrization}
For every $\tau_{\theta,p}\in\Ts$, let
\begin{equation}\label{eq:rejection}
 x:=F(\theta)-p\Delta F(\theta).
\end{equation}
The selected observation under $\tau_{\theta,p}$ has the same distribution
as the selected observation under $\tau_x^F$.
\end{lemma}

\begin{proof}
Take independent uniform variables $U_1,\ldots,U_n$ on $(0,1)$ and write
$X_i=F^{\leftarrow}(U_i)$. Consider the rule that accepts the first
$U_i>x$ with $i<n$, taking $X_n$ if no such index exists.
Since $F(\theta-)\le x\le F(\theta)$, it rejects every $X_i<\theta$
and accepts every $X_i>\theta$, apart from null events. When
$\Delta F(\theta)>0$, the conditional law of $U_i$ given $X_i=\theta$
is uniform on $(F(\theta-),F(\theta)]$, and hence
\[
 \mathbb P(U_i>x\mid X_i=\theta)
 =\frac{F(\theta)-x}{\Delta F(\theta)}=p.
\]
The pairs $(X_i,U_i)$ are independent, so these conditional decisions
have the same joint law as the independent acceptance coins in
$\tau_{\theta,p}$. When the atom has zero mass there is no tie to resolve.
The same construction gives $\tau_x^F$ for $0<x<1$ by \eqref{eq:tie},
and gives its endpoint rules for $x=0,1$. Both rules therefore have the
same selected-value distribution. The argument includes $p=0,1$ and
allows $F$ to be constant on intervals.
\end{proof}

One can also randomize before observing the sequence, by drawing $x$
from a probability measure on $[0,1]$. We call this a mixed quantile
strategy. This is distinct from the independent coins used to resolve
ties after the threshold has been chosen. Although a mixture cannot
improve the best expected reward for a fixed known law, the game below
allows mixtures across quantiles. We shall show in Section~\ref{sec:convex}
that a single quantile can replace any such mixture simultaneously for
all laws $F$.

To compare the two expected rewards, consider their tail probabilities
at a fixed reward level $t\in[0,1]$. Put $y=F(t)$, so that $y$ is the
probability that an observation does not exceed $t$. The rejection
probability of the threshold rule is $x$. If $x>y$, an observation
accepted before the final step necessarily exceeds $t$. If $x\le y$,
an early acceptance need not exceed $t$. These two cases account for
the two branches of the following kernel, identified by
GMZ~\cite[Theorem~2]{GMZ}:
\begin{equation}\label{eq:kernel}
 A_n(x,y)=
 \begin{cases}
 y(x^{n-1}-y^{n-1}),&x>y,\\[1mm]
 (1-y)\displaystyle\sum_{k=1}^{n-1}(y^k-x^k),&x\le y.
 \end{cases}
\end{equation}
Integrating the gap between the tail probabilities of the maximum and
the selected observation gives the regret. The next lemma establishes
this representation.

\begin{lemma}[Tail-integral representation]\label{lem:representation}
For every $F\in\Funit$ and $x\in[0,1]$,
\begin{equation}\label{eq:layercake}
 M_n(F)-V_n(\tau_x^F;F)=\int_0^1 A_n(x,F(t))\,\mathrm dt.
\end{equation}
Equivalently, if $T$ is uniform on $[0,1]$ and $\mu_F$ denotes the
probability law of $F(T)$, then the regret is
$\int A_n(x,y)\,\mathrm d\mu_F(y)$.
\end{lemma}

\begin{proof}
Use the uniform variables from Lemma~\ref{lem:parametrization}. The
selected observation has the same distribution as under the rule that
accepts the first $U_i>x$ before time $n$ and otherwise takes $X_n$.

Fix $t$ and write $y=F(t)$. Apart from null sets,
$\{X_i>t\}=\{U_i>y\}$. Hence the probability that the prophet's observation
exceeds $t$ is $1-y^n$. If $x>y$, every early acceptance exceeds $t$,
and the probability that the selected observation exceeds $t$ is
\[
 (1-x^{n-1})+x^{n-1}(1-y)=1-yx^{n-1}.
\]
If $x\le y$, it is
\[
 (1-y)\sum_{j=0}^{n-2}x^j+x^{n-1}(1-y)
 =(1-y)\sum_{j=0}^{n-1}x^j.
\]
Subtracting these expressions from $1-y^n$ gives \eqref{eq:kernel}.
Integrating the tail probabilities over $t\in[0,1]$ proves
\eqref{eq:layercake}. The endpoint rules satisfy the same identities.
\end{proof}

The law $\mu_F$ may assign mass to $0$ or $1$, for example when the
support of $F$ stays away from an endpoint of $[0,1]$. Such mass contributes
zero to the integral because $A_n(x,0)=A_n(x,1)=0$.

The tail-integral representation separates the rule coordinate $x$
from the law of the observations, which enters only through $\mu_F$.
This is the basis of the difference game of GMZ: the stopping player
chooses a quantile, and the opposing player chooses a probability measure
on the second kernel coordinate. We use $A_n$ as the payoff kernel and
$[0,1]$ as the strategy interval for each player. Let $\lambda$ and $\mu$ be probability measures on
$[0,1]$, and write
\[
 \bar A_n(\lambda,\mu)
 :=\iint A_n(x,y)\,\mathrm d\lambda(x)\,\mathrm d\mu(y).
\]
The minimizing player chooses $\lambda$, a distribution over quantiles.
The maximizing player, called the adversary, chooses $\mu$. A strategy
concentrated at one point is called pure, and $\delta_u$ denotes the
probability measure concentrated at $u$. A pair $(\lambda^*,\mu^*)$ is
a saddle point if
\[
 \bar A_n(\lambda^*,\mu)\le \bar A_n(\lambda^*,\mu^*)
 \le \bar A_n(\lambda,\mu^*)
 \qquad\text{for all }\lambda,\mu.
\]
GMZ~\cite[Theorem~2]{GMZ} prove that
\begin{equation}\label{eq:game}
 v_n=\sup_\mu\inf_\lambda\bar A_n(\lambda,\mu).
\end{equation}
An optimal minimizing strategy bounds the regret from above, and an
optimal maximizing strategy gives the candidate for a matching lower
bound. To establish sharpness for the stopping problem directly, we will
also realize the maximizing measure as $\mu_F$ for a particular law $F$.
We call a distribution least favorable if its smallest attainable regret,
that is, the inner infimum in \eqref{eq:regretdef}, equals $v_n$.

\section{Convexification of the kernel}\label{sec:convex}

To solve the game, it is useful to reduce the minimization over
probability measures to a minimization over individual quantiles. The
probability of reaching the last observation provides a convenient
coordinate: put $z=x^{n-1}$. Written in this coordinate, the first branch
of $A_n$ is affine, and the second is a constant minus a nonnegative
combination of concave powers of $z$. The two expressions meet on the
boundary between the branches, and their maximum is the kernel. This
gives the convexity needed to apply Jensen's inequality.

\begin{lemma}[Convexification]\label{lem:convex}
For fixed $y\in[0,1]$, the function
\[
 B_n(z,y):=A_n(z^{1/(n-1)},y),\qquad 0\le z\le1,
\]
is convex in $z$. Consequently, for any probability measure $\lambda$
on $[0,1]$, the single quantile
\begin{equation}\label{eq:meanquantile}
 \bar x=\left(\int x^{n-1}\,\mathrm d\lambda(x)\right)^{1/(n-1)}
\end{equation}
satisfies
\begin{equation}\label{eq:derandomize}
 A_n(\bar x,y)\le\int A_n(x,y)\,\mathrm d\lambda(x)
 \quad\hbox{for every }y\in[0,1].
\end{equation}
\end{lemma}

\begin{proof}
For all $z,y\in[0,1]$, the two branches can be combined as
\[
 B_n(z,y)=\max\left\{
 yz-y^n,\quad
 (1-y)\sum_{k=1}^{n-1}\bigl(y^k-z^{k/(n-1)}\bigr)
 \right\}.
\]
Indeed, the first expression is nonnegative and the second nonpositive
when $z^{1/(n-1)}\ge y$, and the reverse holds otherwise. The first
expression is affine. The second is convex, because
$z\mapsto z^{k/(n-1)}$ is concave for $1\le k\le n-1$.
The maximum of two convex functions is convex. Jensen's inequality gives
\eqref{eq:derandomize}.
\end{proof}

Integrating \eqref{eq:derandomize} against $\mu_F$ shows that the
replacement works for every underlying distribution $F$.

The single-threshold reduction is a specialization of the construction
of Jiang, Ma and Zhang~\cite[Appendix~E of arXiv:2205.00588v5, equation~(66)]{JMZ}.
When at most one of the first $n-1$ observations may be selected, their
construction preserves the expected number of selections, namely the
average of $1-x^{n-1}$ over the mixed strategy. The replacement therefore has
quantile $\bar x$ as in
\eqref{eq:meanquantile}. The contribution from the compulsory last
observation is $x^{n-1}\mathbb E_F X_n$, which is preserved as well.
For $n=2$ the convexity argument above is elementary. We use this direct
kernel proof because the same convex coordinate will establish the
matching lower bound in Section~\ref{sec:saddle}. Loiseau, Molina, Perchet, Perez-Salazar and
Verdugo~\cite{LMPPV} study related quantile kernels for ratio-type prophet
inequalities under additional variance constraints; their objective and
stopping convention are different.

\section{Maximizing the two branches}

Lemma~\ref{lem:convex} reduces the problem to minimizing the largest
payoff over $y$ for a single quantile $x$. Recall that $x$ is the rejection
probability and $y=F(t)$ for a reward level $t$. In the range $y<x$,
an early acceptance necessarily exceeds $t$, but the rule may be forced
to take a smaller final observation. In the range $y\ge x$, it may
also accept a value no greater than $t$ before the final observation.
The largest payoff in each range can be computed by differentiation
with respect to $y$. Denote these maxima by $D_n(x)$ and $U_n(x)$,
respectively. The constants needed for the first branch are
\begin{equation}\label{eq:alphabeta}
 \beta_n:=n^{-1/(n-1)},\qquad
 \alpha_n:=\left(1-\frac1n\right)\beta_n.
\end{equation}
For the second branch, the derivative involves the average of the first
$n$ powers of $x$, including the zeroth power. For $x\in[0,1]$, put
\begin{equation}\label{eq:mz}
 m_n(x):=\frac1n\sum_{j=0}^{n-1}x^j,
 \qquad y_-(x):=\beta_nx,
 \qquad y_+(x):=m_n(x)^{1/(n-1)}.
\end{equation}
Thus $m_n(1)=y_+(1)=1$, and
$m_n(x)=(1-x^n)/(n(1-x))$ when $x<1$.
The points $y_-(x)$ and $y_+(x)$ defined in \eqref{eq:mz} are the
candidate maximizers in the two ranges. Their payoffs are
\begin{align}
 D_n(x)&:=\alpha_nx^n,\label{eq:D}\\
 U_n(x)&:=1-m_n(x)\bigl[n-(n-1)y_+(x)\bigr].\label{eq:U}
\end{align}

\begin{lemma}[Exact branch maximizers]\label{lem:branches}
For every $x\in(0,1)$,
\begin{align}
 \max_{0\le y<x} A_n(x,y)&=D_n(x),\label{eq:lowerbranch}\\
 \max_{x\le y\le1} A_n(x,y)&=U_n(x).\label{eq:upperbranch}
\end{align}
The unique maximizing points are $y_-(x)$ and $y_+(x)$, respectively,
and satisfy $y_-(x)<x<y_+(x)$.
\end{lemma}

\begin{proof}
On the interval $y<x$ the payoff is $yx^{n-1}-y^n$, with derivative
$x^{n-1}-ny^{n-1}$. The derivative changes from positive to negative at
$y=\beta_nx$, so this point is the unique maximum. Substitution gives
\[
 A_n(x,\beta_nx)=\beta_nx^n(1-1/n)=D_n(x).
\]
For $y\ge x$, summing the geometric progression in the kernel gives
\[
 A_n(x,y)=y-y^n-(1-y)\sum_{k=1}^{n-1}x^k.
\]
Its derivative is therefore
\begin{equation}\label{eq:dyupper}
 \partial_y A_n(x,y)
 =1+\sum_{k=1}^{n-1}x^k-ny^{n-1}
 =n\bigl[m_n(x)-y^{n-1}\bigr].
\end{equation}
For $0<x<1$, the average $m_n(x)$ lies strictly between $x^{n-1}$
and $1$. Thus $x<y_+(x)<1$, and the derivative changes sign there.
Using $y_+(x)^{n-1}=m_n(x)$ in the preceding expression gives $U_n(x)$.
\end{proof}

Increasing the rejection probability raises the first maximum and lowers
the second. The next lemma establishes that the two maxima meet exactly
once, so their balance is characterized by one scalar equation.

\begin{lemma}[Monotone crossing]\label{lem:crossing}
$D_n$ is strictly increasing on $(0,1)$ and $U_n$ is strictly decreasing on $(0,1)$.  In addition,
\[
 D_n(0)=0,\qquad U_n(0)=\alpha_n,
 \qquad
 D_n(1)=\alpha_n,\qquad U_n(1)=0.
\]
Consequently there is a unique $x_n^*\in(0,1)$ such that
\begin{equation}\label{eq:scalar}
  \alpha_n (x_n^*)^n
 =1-m_n(x_n^*)
 \Bigl[n-(n-1)m_n(x_n^*)^{1/(n-1)}\Bigr].
\end{equation}

Denote the common value at the crossing by
\begin{equation}\label{eq:qdef}
 q_n:=D_n(x_n^*)=U_n(x_n^*)
 =\alpha_n(x_n^*)^n.
\end{equation}
\end{lemma}

\begin{proof}
The derivative of $D_n(x)=\alpha_nx^n$ is positive. For $U_n$, apply
the chain rule to $A_n(x,y_+(x))$. The term involving $y_+'(x)$ vanishes
because $\partial_y A_n(x,y_+(x))=0$, leaving
\begin{equation}\label{eq:Uderivative}
 U_n'(x)
 =-(1-y_+(x))\sum_{k=1}^{n-1}k x^{k-1}<0.
\end{equation}
Evaluation at $0$ and $1$ gives the stated endpoint values. Thus
$D_n-U_n$ is continuous, strictly increasing, and changes sign once.
\end{proof}

Thus a deterministic quantile $x$ has worst kernel regret
\begin{equation}\label{eq:pureenv}
 \sup_{0\le y\le1}A_n(x,y)=\max\{D_n(x),U_n(x)\}.
\end{equation}
Below $x_n^*$ the larger loss is $U_n(x)$, and raising $x$ reduces it.
Above $x_n^*$ the larger loss is $D_n(x)$, and lowering $x$ reduces it.
Their intersection therefore minimizes the worst loss. Mixtures cannot
improve it by Lemma~\ref{lem:convex}. To obtain a sharp bound for the original stopping problem, however, the
upper estimate must be matched by a single adversarial measure against
all quantiles. The next section constructs such a measure.

\section{The finite-horizon saddle point}\label{sec:saddle}

At the crossing, both branch maximizers yield the same payoff. An
adversary who uses just one of these points would allow the stopping
player to reduce the payoff by moving the quantile away from the
crossing. We instead put mass at both points, choosing the weights so
that the derivative of the average payoff vanishes there. Convexity in
$x^{n-1}$ then ensures a global minimum.

Denote the two maximizing points at the crossing by
\begin{equation}\label{eq:ystars}
 y_{-,n}^*:=\beta_n x_n^*,
 \qquad
 y_{+,n}^*:=y_+(x_n^*).
\end{equation}
Their payoff derivatives, with respect to $x$ at $x_n^*$, have
opposite signs. Denote their positive magnitudes by $P_n$ and $Q_n$;
let $\rho_n$ be the probability assigned to the upper point:
\begin{align}
 P_n&:=(n-1)y_{-,n}^*(x_n^*)^{n-2},\label{eq:P}\\
 Q_n&:=(1-y_{+,n}^*)\sum_{k=1}^{n-1}k(x_n^*)^{k-1},\label{eq:Q}\\
 \rho_n&:=\frac{P_n}{P_n+Q_n}.\label{eq:rho}
\end{align}
Here $P_n$ and $-Q_n$ are the slopes of
$x\mapsto A_n(x,y_{-,n}^*)$ and $x\mapsto A_n(x,y_{+,n}^*)$ at
$x_n^*$, respectively. Both $P_n$ and $Q_n$ are positive, and $\rho_n$ is
exactly the weight for which $(1-\rho_n)P_n-\rho_n Q_n=0$.
The payoff against the proposed adversarial strategy is
\begin{equation}\label{eq:Gdef}
 G_n(x):=\rho_n A_n(x,y_{+,n}^*)
       +(1-\rho_n)A_n(x,y_{-,n}^*).
\end{equation}

\begin{theorem}[Exact finite-horizon solution]\label{thm:exact}
For every $n\ge2$,
\begin{equation}\label{eq:exactvalue}
  v_n=q_n
 =\alpha_n(x_n^*)^n,
\end{equation}
where $x_n^*$ is the unique solution of \eqref{eq:scalar}.

An optimal decision-maker strategy in the game \eqref{eq:game} is pure:
\[
 \lambda_n^*=\delta_{x_n^*}.
\]
An optimal strategy for the adversary is the two-point measure
\begin{equation}\label{eq:mustar}
  \mu_n^*
 =\rho_n\delta_{y_{+,n}^*}
 +(1-\rho_n)\delta_{y_{-,n}^*}.
\end{equation}
Thus $(\lambda_n^*,\mu_n^*)$ is a saddle point of the difference game.
\end{theorem}

\begin{proof}
Lemmas~\ref{lem:branches}--\ref{lem:crossing} give
\begin{equation}\label{eq:uppergame}
 A_n(x_n^*,y)\le q_n \quad\hbox{for all }y\in[0,1].
\end{equation}
For the lower bound, $G_n$ from \eqref{eq:Gdef} is differentiable at
$x_n^*$ because $y_{-,n}^*<x_n^*<y_{+,n}^*$. Its derivative is
\[
 G_n'(x_n^*)=(1-\rho_n)P_n-\rho_n Q_n=0.
\]
By Lemma~\ref{lem:convex}, the function
$z\mapsto G_n(z^{1/(n-1)})$ is convex. Its derivative at the interior
point $z=(x_n^*)^{n-1}$ is zero, so this point is a global minimizer.
Both support points attain the branch maximum $q_n$ at $x_n^*$;
therefore
\begin{equation}\label{eq:lowergame}
 G_n(x)\ge G_n(x_n^*)=q_n \quad\hbox{for every }x\in[0,1].
\end{equation}
The two inequalities establish the saddle property, also against mixed
strategies by integration. For the original stopping problem,
\eqref{eq:uppergame} and Lemma~\ref{lem:representation} give the bound
$v_n\le q_n$. For the reverse inequality, take a distribution with probabilities
$y_{-,n}^*$, $y_{+,n}^*-y_{-,n}^*$, and $1-y_{+,n}^*$ at
$0$, $1-\rho_n$, and $1$, respectively. Its cumulative distribution
function equals $y_{-,n}^*$ on an interval of length $1-\rho_n$ and
$y_{+,n}^*$ on an interval of length $\rho_n$. Its regret against any
quantile $x$ is therefore $G_n(x)\ge q_n$. Every rule in $\Ts$ is
equivalent to a quantile rule by Lemma~\ref{lem:parametrization}, so
$v_n\ge q_n$.
\end{proof}

\begin{remark}[Computational consequence]
Bisection applied to $D_n-U_n$ computes the sharp constant: this function
is strictly increasing and has opposite signs at $0$ and $1$. No matrix
game needs to be discretized.
\end{remark}

\begin{remark}[Uniqueness]\label{rem:uniqueness}
For $n\ge3$, the minimizing game strategy is uniquely
$\delta_{x_n^*}$. To see this, set $H_n(z)=G_n(z^{1/(n-1)})$ and
$z_*=(x_n^*)^{n-1}$. This is a convex function with a minimum at $z_*$.
Near $z_*$, the upper-point contribution contains the term
$-\rho_n(1-y_{+,n}^*)z^{1/(n-1)}$, which is strictly convex; its other
terms are convex. Hence $H_n$ is strictly convex near $z_*$. A second
minimum would force it to be constant on the segment joining the two
minima, contradicting this local strict convexity. Thus
$G_n(x)>q_n$ whenever $x\ne x_n^*$.

If $\lambda\ne\delta_{x_n^*}$, it puts positive mass outside $x_n^*$,
so
\begin{equation}\label{eq:strictmixed}
 \sup_\mu\bar A_n(\lambda,\mu)
 \ge\bar A_n(\lambda,\mu_n^*)
 =q_n+\int_{[0,1]}(G_n(x)-q_n)\,\mathrm d\lambda(x)>q_n.
\end{equation}
The last integral is strictly positive because its integrand is
nonnegative and positive off $x_n^*$. This proves uniqueness among all
mixed quantile strategies.

For $n=2$, the unique pure optimum is $x_2^*=1/2$, but mixed optima
also exist. For example, $(\delta_{2/5}+\delta_{3/5})/2$ has payoff
$y(1/2-y)$ for $y\le2/5$, $(1-y)(y-1/2)$ for $y\ge3/5$, and
$y-y^2-1/5\le1/20$ in between. Its maximum is $1/16$, attained at
$y=1/4$ and $y=3/4$.
\end{remark}

\section{Optimal rule and least-favorable distribution}\label{sec:stopping}

Theorem~\ref{thm:exact} identifies a quantile level $x_n^*$ that
attains the minimax guarantee uniformly over all laws. The numerical
threshold still depends on $F$, and the rule optimized for a particular
law may have a smaller regret than this universal rule. We now state the guarantee in terms of the original
observations and describe the law that makes it sharp.

The two-point adversarial measure also has a direct interpretation
in the observation model. Since $\mu_F$ is the law of $F(T)$ for a
uniform variable $T$, its masses are the lengths of the intervals on
which the distribution function is constant. Two such intervals lead
to a law supported on three points: the weights of $\mu_n^*$ are the
lengths of the two constant pieces of the distribution function, while
the probabilities of the three observation values are its jump sizes.

\begin{corollary}[Exact universal rule]\label{cor:rule}
For every $F\in\Funit$, choose
\begin{equation}\label{eq:univtheta}
 \theta_n(F):=F^{\leftarrow}(x_n^*)
\end{equation}
and use the tie probability \eqref{eq:tie} with $x=x_n^*$.  Then
\begin{equation}\label{eq:univbound}
 M_n(F)-V_n(\tau_{x_n^*}^F;F)\le q_n.
\end{equation}
The constant $q_n$ is best possible over randomized single-threshold rules.
\end{corollary}

\begin{proof}
By \eqref{eq:layercake} and \eqref{eq:uppergame},
\[
 M_n(F)-V_n(\tau_{x_n^*}^F;F)
 =\int_0^1 A_n(x_n^*,F(t))\,\mathrm dt\le q_n.
\]
Sharpness is witnessed by the distribution constructed in the proof of
Theorem~\ref{thm:exact}.
\end{proof}

\begin{corollary}[Explicit three-point least-favorable distribution]\label{cor:lfd}
Let $X$ have distribution $F_n^*$ supported on
\begin{equation}\label{eq:threepoints}
 \{0,\,1-\rho_n,\,1\}
\end{equation}
with probabilities
\begin{align}
 \mathbb P_{F_n^*}(X=0)&=y_{-,n}^*,\label{eq:p0}\\
 \mathbb P_{F_n^*}(X=1-\rho_n)&=y_{+,n}^*-y_{-,n}^*,\label{eq:pmid}\\
 \mathbb P_{F_n^*}(X=1)&=1-y_{+,n}^*.\label{eq:p1}
\end{align}
Then $F_n^*$ is least favorable for the difference regret.  Its quantile function is
\begin{equation}\label{eq:qfun}
 (F_n^*)^{\leftarrow}(u)=
 \begin{cases}
 0,&0\le u\le y_{-,n}^*,\\
 1-\rho_n,&y_{-,n}^*<u\le y_{+,n}^*,\\
 1,&y_{+,n}^*<u\le1,
 \end{cases}
\end{equation}
with the endpoint values assigned as displayed and the values at the
interior jumps chosen according to the left-continuous inverse convention.

The quantile $x_n^*$ lies between $y_{-,n}^*$ and $y_{+,n}^*$, so an
optimal rule for this law uses threshold $1-\rho_n$ and acceptance
probability
\begin{equation}\label{eq:ptieexact}
 p_n^*=
 \frac{y_{+,n}^*-x_n^*}{y_{+,n}^*-y_{-,n}^*}.
\end{equation}
\end{corollary}

\begin{proof}
For $0\le t<1-\rho_n$ the distribution function is $y_{-,n}^*$;
for $1-\rho_n\le t<1$ it is $y_{+,n}^*$.
Thus the law of $F_n^*(T)$, with $T$ uniform, is exactly $\mu_n^*$.
Lemma~\ref{lem:representation} now gives
\[
 M_n(F_n^*)-V_n(\tau_x^{F_n^*};F_n^*)=G_n(x)\ge q_n
\]
for every $x$, with equality at $x_n^*$. Lemma~\ref{lem:parametrization}
covers every rule in $\Ts$, proving that $F_n^*$ is least favorable.
The formula for $p_n^*$ follows from \eqref{eq:tie} because the atom at
$1-\rho_n$ has size $y_{+,n}^*-y_{-,n}^*$.
\end{proof}

The construction above establishes existence of a least-favorable law.
The equality conditions in the branch maxima also determine whether
other laws can have the same smallest attainable regret.

\begin{corollary}[Uniqueness of the least-favorable distribution]\label{cor:unique_lfd}
For every $n\ge2$, the measure $\mu_n^*$ is the unique optimal strategy
of the maximizing player in the difference game. Moreover, $F_n^*$ is
the unique distribution in $\Funit$ satisfying
\[
 \inf_{\tau\in\Ts}
 \{M_n(F)-V_n(\tau;F)\}=q_n.
\]
\end{corollary}

\begin{proof}
Fix $n\ge2$ and write $x^*=x_n^*$, $a=y_{-,n}^*$, $b=y_{+,n}^*$,
and $q=q_n$. By Lemma~\ref{lem:branches} and the defining equality of the branch
maxima, $A_n(x^*,y)\le q$, with equality exactly at $y=a,b$.
If $\mu$ is an optimal maximizing strategy, then
\[
 q\le\int A_n(x^*,y)\,\mathrm d\mu(y)\le q.
\]
Consequently, $\mu=(1-w)\delta_a+w\delta_b$ for some $w\in[0,1]$.
The function $x\mapsto\int A_n(x,y)\,\mathrm d\mu(y)$ has a minimum at
$x^*\in(a,b)$ and is differentiable there. Thus
\[
 0=(1-w)P_n-wQ_n,
\]
which gives $w=\rho_n$ and proves uniqueness of $\mu_n^*$.

If $F$ is least favorable, Lemmas~\ref{lem:parametrization} and
\ref{lem:representation} give $\int A_n(x,y)\,\mathrm d\mu_F(y)\ge q_n$
for every $x$. Integrating with respect to any $\lambda$ shows that
$\mu_F$ is an optimal maximizing strategy, so $\mu_F=\mu_n^*$.
Hence $F(t)$ takes the values $a$ and $b$ for almost every
$t\in(0,1)$, on sets of lengths $1-\rho_n$ and $\rho_n$,
respectively. Monotonicity and right continuity of $F$ then give
\[
 F(t)=
 \begin{cases}
 a,&0\le t<1-\rho_n,\\
 b,&1-\rho_n\le t<1,
 \end{cases}
 \qquad F(1)=1.
\]
This is precisely the distribution function of $F_n^*$.
\end{proof}

\begin{remark}
Uniqueness here concerns the least-favorable law after optimizing the
rule for that law. The laws at which the fixed universal quantile
$x_n^*$ has regret $q_n$ need not be unique. In fact, for any
$d\in[0,1]$, placing probabilities $y_{-,n}^*$,
$y_{+,n}^*-y_{-,n}^*$, and $1-y_{+,n}^*$ at $0,d,1$ gives regret
$q_n$ at $x_n^*$. Only $d=1-\rho_n$ gives a least-favorable law.
\end{remark}

\begin{remark}[Two observations]
When $n=2$, the formulas reduce to
\[
 x_2^*=\tfrac12,\quad y_{-,2}^*=\tfrac14,\quad
 y_{+,2}^*=\tfrac34,\quad \rho_2=\tfrac12,\quad q_2=\tfrac1{16}.
\]
The least-favorable law assigns probabilities $1/4,1/2,1/4$ to
$0,1/2,1$ and has mean $1/2$. With two observations, it is best to
reject an initial $0$, accept an initial $1$, and do either at $1/2$:
accepting and continuing both have conditional value $1/2$. This gives
a stopping reward of $5/8$, whereas
$M_2=1-(3/4)^2+\tfrac12((3/4)^2-(1/4)^2)=11/16$.
Their difference is $1/16$. Any tie probability at $1/2$ is optimal
for this particular law; the displayed $p_2^*=1/2$ is one choice.
\end{remark}

\section{Large-horizon asymptotics}\label{sec:asympt}

The formula for $q_n$ contains the factor $(x_n^*)^n$. For any fixed
$x\in[0,1)$ the corresponding lower-branch maximum tends to zero, so an
appropriate scale for investigating a nonzero limit is $x=e^{-s/n}$
with $s>0$. In this coordinate the probability of reaching the last
observation is $e^{-s(n-1)/n}$, which remains nondegenerate as the horizon
grows with $s$ fixed. We will show that the optimal rescaled parameter
converges to a positive constant, and use the balance of the two branches
to calculate its first correction. All logarithms are natural. Write
\begin{equation}\label{eq:sndef}
 s_n:=-n\log x_n^*.
\end{equation}
The limiting crossing derived in the proof below is determined by the
nonzero root $r_*>1$ of
\begin{equation}\label{eq:rstar}
 e^{r_*}=er_*+1.
\end{equation}
This solution is unique: $e^r-er-1$ decreases on $(0,1)$, increases on
$(1,\infty)$, is negative at $1$, and tends to infinity. The equation also
has the root $r=0$; it does not solve the limiting crossing equation and plays
no role below (see the proof of Theorem~\ref{thm:asympt}). The constants in the expansion can all be expressed in terms of this
root. Write
\begin{align}
 a_*&:=r_*-1, & c_*&:=e^{-r_*},\label{eq:acstar}\\
 \rho_*&:=\frac{r_*}{1+e(r_*-1)}, &
 \gamma_*&:=1-\rho_*+\frac{\rho_*ea_*}{2}.\label{eq:gammastar}
\end{align}
Since $a_*>0$ and $1+ea_*=r_*+(e-1)a_*$, we have $0<\rho_*<1$.
Numerically,
\[
 \begin{aligned}
 r_*&=1.750786722680\ldots,& c_*&=0.173637285311\ldots,\\
 \rho_*&=0.575755718697\ldots,& \gamma_*&=1.011759783295\ldots.
 \end{aligned}
\]

\begin{theorem}[Limit and second-order expansion]\label{thm:asympt}
As $n\to\infty$,
\begin{equation}\label{eq:snexp}
 -n\log x_n^*
 =r_*
 -\frac{\rho_*\log n+1-\gamma_*}{n}
 +O\!\left(\frac{(\log n)^2}{n^2}\right).
\end{equation}
The sharp regret satisfies
\begin{equation}\label{eq:qexp}
  q_n
 =c_*\left[
 1-\frac{(1-\rho_*)\log n+\gamma_*}{n}
 +O\!\left(\frac{(\log n)^2}{n^2}\right)
 \right].
\end{equation}
In particular,
\begin{equation}\label{eq:limit}
  \lim_{n\to\infty}
 v_n
 =c_*=0.17363728531099157\ldots.
\end{equation}
\end{theorem}

\begin{proof}
We first identify the limiting crossing, and then expand the finite
crossing about it. For $s>0$, put $x=e^{-s/n}$. The lower-branch maximum is
\begin{equation}\label{eq:Dscaled}
 D_n(e^{-s/n})=\alpha_n e^{-s},
 \qquad
 \alpha_n
 =1-\frac{\log n+1}{n}
 +O\!\left(\frac{(\log n)^2}{n^2}\right).
\end{equation}
The upper branch has its own maximizing point. To compare it with the
rescaled quantile, write that point as
\[
 y_+(e^{-s/n})=e^{-t_{+,n}(s)/n}.
\]
From $y_+(e^{-s/n})^{n-1}=m_n(e^{-s/n})$,
\begin{equation}\label{eq:tnexact}
 t_{+,n}(s)
 =\frac{n}{n-1}
 \log\frac{n(1-e^{-s/n})}{1-e^{-s}}.
\end{equation}
Define
\begin{equation}\label{eq:t0}
 h(s):=\frac{1-e^{-s}}s,\qquad t_0(s):=-\log h(s).
\end{equation}
Equation~\eqref{eq:tnexact} gives
$t_{+,n}(s)=t_0(s)+O(n^{-1})$, uniformly on compact subintervals of
$(0,\infty)$. Since $e^{-s}<h(s)<1$, we have $0<t_0(s)<s$.
For $s>0$ and $0\le t\le s$, direct summation of the kernel gives the
exact identity
\begin{equation}\label{eq:scaledexact}
 A_n(e^{-s/n},e^{-t/n})
 =1-e^{-t}-(1-e^{-s})\frac{1-e^{-t/n}}{1-e^{-s/n}}.
\end{equation}
Together with the first branch of \eqref{eq:kernel}, this identifies
the limit, for fixed $s,t\ge0$, as
\begin{equation}\label{eq:K}
 K(s,t)=
 \begin{cases}
 e^{-s}-e^{-t},&s<t,\\[1mm]
 1-e^{-t}-\dfrac ts(1-e^{-s}),&s\ge t.
 \end{cases}
\end{equation}
We set $K(0,0)=0$, the continuous extension at the origin.
For $s>0$, the function $t\mapsto K(s,t)$ on $[0,s]$ has its unique
maximum at $t_0(s)$, because
\[
 \partial_t K(s,t)=e^{-t}-\frac{1-e^{-s}}s.
\]
Hence its limiting envelope is
\begin{equation}\label{eq:Ulim}
 U(s)=1-e^{-t_0(s)}
 -\frac{t_0(s)}s(1-e^{-s}).
\end{equation}
The limiting crossing equation $U(s)=e^{-s}$ simplifies to
\[
 s=1+t_0(s),
\]
which is equivalent to $e^s=es+1$ for $s>0$. The solution $s=0$ of the latter algebraic equation is extraneous to the crossing, since
$\lim_{s\downarrow0}(U(s)-e^{-s})=-1$.
To pass from the limiting equation to the finite roots, note that
$h'(s)<0$ for $s>0$ and
\[
 U(s)=1-h(s)+h(s)\log h(s),\qquad U'(s)=h'(s)\log h(s)>0.
\]
Therefore $F_0(s):=U(s)-e^{-s}$ is strictly increasing, with limits
$-1$ and $1$ at $0$ and $\infty$. It has the unique zero $r_*$.
The finite functions $U_n(e^{-s/n})-D_n(e^{-s/n})$ are also strictly
increasing in $s$. For any small fixed $\varepsilon>0$, pointwise
convergence at $r_*-\varepsilon$ and $r_*+\varepsilon$ gives opposite
signs for all sufficiently large $n$. Their unique roots are therefore
trapped between these two points. Thus
$-n\log x_n^*\to r_*$ and $q_n\to c_*$.

It remains to calculate the displacement of the finite crossing from
$r_*$. For this we need uniform expansions in a neighborhood of $r_*$.
Choose $\eta>0$ small enough that the interval
$I=[r_*-\eta,r_*+\eta]$ is contained in $(0,\infty)$ and $F_0'$ is
bounded below by a positive constant on $I$. The root argument puts
$s_n$ in $I$ for all sufficiently large $n$. Uniformly for $s\in I$
and $0\le t\le s$,
\begin{equation}\label{eq:ratioexp}
 \frac{1-e^{-t/n}}{1-e^{-s/n}}
 =\frac ts\left[1+\frac{s-t}{2n}+O(n^{-2})\right].
\end{equation}
Since $I$ stays away from zero, the normalized denominator
$(1-e^{-s/n})/(s/n)$ stays bounded away from zero for large $n$;
Taylor's formula on the bounded range $0\le t\le r_*+\eta$ gives
a common constant in the remainder. Define
\begin{equation}\label{eq:u1}
 u_1(s):=-\frac{(1-e^{-s})t_0(s)[s-t_0(s)]}{2s}.
\end{equation}
Substituting \eqref{eq:ratioexp} into \eqref{eq:scaledexact} and evaluating
at $t=t_{+,n}(s)$ gives
\begin{equation}\label{eq:Unexp}
 U_n(e^{-s/n})=U(s)+\frac{u_1(s)}n+O(n^{-2}).
\end{equation}
Indeed, $t_{+,n}(s)-t_0(s)=O(n^{-1})$ and
$\partial_tK(s,t_0(s))=0$, so shifting the maximizing point contributes
only $O(n^{-2})$. Replacing $t_{+,n}(s)$ by $t_0(s)$ in the coefficient
of $1/n$ has an error of the same order.
At $s=r_*$, $t_0(r_*)=a_*=r_*-1$ and
\begin{equation}\label{eq:u1r}
 u_1(r_*)=-\frac{e a_* c_*}{2}.
\end{equation}

For the derivative of $F_0(s)=U(s)-e^{-s}$ at its root, use
$h(r_*)=ec_*$ and $h'(r_*)=-(e-1)c_*/r_*$. Then
\begin{equation}\label{eq:Fprime}
 F_0'(r_*)=c_*+\frac{a_*(e-1)c_*}{r_*}
 =\frac{c_*}{\rho_*}>0.
\end{equation}
The expansions on $I$ imply
\[
 \sup_{s\in I}\left|U_n(e^{-s/n})-D_n(e^{-s/n})-F_0(s)\right|
 =O((\log n)/n).
\]
The mean value theorem first gives $s_n-r_*=O((\log n)/n)$.
The quadratic Taylor remainder is then $O((\log n)^2/n^2)$;
$u_1$ and the other coefficient functions have bounded derivatives on
$I$.
At the exact crossing, \eqref{eq:Dscaled} and \eqref{eq:Unexp} give
\[
 0=F_0(s_n)+\frac{u_1(s_n)+e^{-s_n}(\log n+1)}n
   +O\!\left(\frac{(\log n)^2}{n^2}\right).
\]
Expanding at $r_*$, using \eqref{eq:u1r} and \eqref{eq:Fprime}, yields
\[
 s_n-r_*
 =-\frac{\rho_*}{n}
 \left(\log n+1-\frac{e a_*}{2}\right)
 +O\!\left(\frac{(\log n)^2}{n^2}\right),
\]
which is \eqref{eq:snexp}.  Finally, substituting this into the exact identity
$q_n=\alpha_ne^{-s_n}$ and using \eqref{eq:Dscaled} gives \eqref{eq:qexp}.
\end{proof}

The asymptotics of the least-favorable law require both maximizing
points, not just the optimal quantile. Their distances from the endpoint
$1$ have different orders, as the following proposition shows.

\begin{proposition}[Asymptotics of the least-favorable law]\label{prop:geometry}
The two support points of the adversary and the mixing weight satisfy
\begin{align}
 -n\log y_{+,n}^*&\longrightarrow a_*=r_*-1,\label{eq:yplusas}\\
 -n\log y_{-,n}^*
 &= -n\log x_n^*+\frac{n}{n-1}\log n
 =\log n+r_*+o(1),\label{eq:yminusas}\\
 \rho_n&\longrightarrow\rho_*.
\label{eq:rhoas}
\end{align}
Equivalently, $1-y_{+,n}^*$ is of order $1/n$, whereas
$1-y_{-,n}^*$ is of order $(\log n)/n$. These are the two boundary
scales at the endpoint $1$.
\end{proposition}

\begin{proof}
Equation~\eqref{eq:tnexact} gives the first limit, while
$y_{-,n}^*=n^{-1/(n-1)}x_n^*$ gives the second after taking logarithms
and using \eqref{eq:snexp}. For the mixing weight, the Riemann sum
\[
 n^{-2}\sum_{k=1}^{n-1}k(x_n^*)^{k-1}
 \longrightarrow \int_0^1 u e^{-r_*u}\,\mathrm du
 =\frac{1-(1+r_*)c_*}{r_*^2}
\]
together with $n(1-y_{+,n}^*)\to a_*$ gives
\[
 P_n/n\to c_*,\qquad Q_n/n\to
 \frac{a_*[1-(1+r_*)c_*]}{r_*^2}
 =\frac{a_*c_*(e-1)}{r_*}.
\]
Taking the ratio $P_n/(P_n+Q_n)$ yields \eqref{eq:rhoas}.
\end{proof}

\section{Monotonicity and a uniform bound}\label{sec:uniform}

The asymptotic expansion gives the direction from which $q_n$ approaches
$c_*$ for sufficiently large $n$, but does not by itself establish
monotonicity for every horizon. Both that question and the search for
an explicit rule that avoids solving the scalar equation can be addressed
by comparing the rescaled kernels. We first give a bound valid for all
values of their two arguments.

\begin{lemma}[Global kernel domination]\label{lem:domination}
For every $n\ge2$ and $s,t\ge0$,
\begin{equation}\label{eq:domination}
 A_n(e^{-s/n},e^{-t/n})\le K(s,t),
\end{equation}
where $K$ is given in \eqref{eq:K}, with $K(0,0)=0$ by continuity.
\end{lemma}

\begin{proof}
Put $x=e^{-s/n}$ and $y=e^{-t/n}$. At $s=t=0$, both kernels are zero.
Otherwise, if $s<t$, then $x>y$ and
\[
 A_n=e^{-s}e^{(s-t)/n}-e^{-t}
 \le e^{-s}-e^{-t}=K(s,t).
\]
If $s\ge t$, direct algebra gives
\[
 A_n-K
 =(1-e^{-s})
 \left[
 \frac ts-\frac{1-e^{-t/n}}{1-e^{-s/n}}
 \right]\le0,
\]
because $u\mapsto(1-e^{-u})/u$ is decreasing on $(0,\infty)$.
\end{proof}

The same rescaling permits comparison of successive horizons. Strict
increase of the kernel away from its zeros will imply strict increase of
the optimized constants, by evaluating at the optimum for the larger
horizon.

\begin{proposition}[Strict increase with the horizon]\label{prop:monotonicity}
The sharp constants in Theorem~\ref{thm:exact} satisfy
\[
 q_2<q_3<\cdots<c_*.
\]
\end{proposition}

\begin{proof}
For $s,t\ge0$, put
\[
 K_n(s,t)=A_n(e^{-s/n},e^{-t/n}),
 \qquad W_n(s)=\sup_{t\ge0}K_n(s,t).
\]
For each $s>0$, the supremum is a positive maximum attained at a
finite point on one of the two branches. Neither $t=0$ nor $t=s$
is a maximizer, since the kernel vanishes there. The endpoint $y=0$
also contributes zero. The finite-horizon solution gives
$q_n=\min_{s>0}W_n(s)$, attained at $s_n=-n\log x_n^*$.

If $s<t$, then
\[
 K_n(s,t)=e^{-s}e^{(s-t)/n}-e^{-t}
\]
is strictly increasing in $n$. If $0<t<s$, then
\[
 K_n(s,t)=1-e^{-t}-(1-e^{-s})
 \frac{1-e^{-t/n}}{1-e^{-s/n}}.
\]
For $R(u)=(1-e^{-tu})/(1-e^{-su})$, $u>0$, we have
\[
 \frac{d}{du}\log R(u)
 =\frac{t}{e^{tu}-1}-\frac{s}{e^{su}-1}>0,
\]
because $v/(e^v-1)$ is strictly decreasing for $v>0$ (the numerator of
its derivative is $e^v(1-v)-1<0$).
Thus the ratio in the preceding display decreases strictly with $n$,
and $K_n(s,t)$ increases strictly on this branch as well.

Fix $s=s_{n+1}$ and let $t_n$ maximize $K_n(s_{n+1},t)$. By the
properties above, strict comparison applies at $t_n$. Therefore
\[
 \begin{aligned}
 q_{n+1}
 &=W_{n+1}(s_{n+1})\\
 &\ge K_{n+1}(s_{n+1},t_n)\\
 &>K_n(s_{n+1},t_n)\\
 &=W_n(s_{n+1})\ge q_n.
 \end{aligned}
\]
Finally, Theorem~\ref{thm:asympt} gives $q_n\to c_*$, so strict increase implies
$q_n<c_*$ for every finite $n$.
\end{proof}

The limiting crossing also supplies a simple choice of quantile for
every finite horizon, with a bound independent of $n$.

\begin{corollary}[Explicit all-$n$ rule]
For every $n\ge2$, choosing the quantile
\begin{equation}\label{eq:simplex}
 \widetilde x_n=e^{-r_*/n}
\end{equation}
guarantees worst-case additive regret at most $c_*$.
\end{corollary}

\begin{proof}
By Lemma~\ref{lem:domination},
\[
 \sup_{0<y\le1} A_n(e^{-r_*/n},y)
 \le \sup_{t\ge0}K(r_*,t).
\]
For $t>r_*$, $K(r_*,t)=c_*-e^{-t}<c_*$. On $[0,r_*]$, its
maximum occurs at $t=a_*$ and equals $c_*$. The endpoint $y=0$ also
has zero payoff, which completes the bound.
\end{proof}

\section{Numerical values and computation}

Bisection of the scalar crossing equation needs only the two branch
maxima. For large $n$ it is convenient to solve for $s=-n\log x$,
rather than for a quantile increasingly close to $1$. To evaluate the
maxima near their crossing, write $\operatorname{expm1}(u)=e^u-1$ for
the elementary function evaluated directly for small arguments. Use
\[
 m_n(e^{-s/n})=\frac{-\operatorname{expm1}(-s)}{-n\operatorname{expm1}(-s/n)}
 \qquad(s>0),
\]
and, writing $m=m_n(e^{-s/n})$ and $\ell=\log(m)/(n-1)$, evaluate the
upper envelope as
\[
 U_n(e^{-s/n})=1-m\{1-(n-1)\operatorname{expm1}(\ell)\}.
\]
At $s=0$ one uses $m_n(1)=1$ and $U_n(1)=0$. Very close to zero,
the expression for $U_n$ may still suffer cancellation in fixed-precision
arithmetic; increased precision or a local series expansion is then
needed if relative accuracy is required.

Table~\ref{tab:numerics} reports $q_n=v_n$ for the difference game on the
full interval $[0,1]$. The entries were computed from \eqref{eq:scalar}
with 100-digit arithmetic and rounded for display. For comparison, at
$n=25$ the discretized value and upper bound reported in Table~1 of
\cite{GMZ} are $0.1572$ and $0.158$; at $n=3000$ they are $0.1731$
and $0.288$. The scalar equation gives the sharp constants in both
cases without increasing the dimension of a discretized game.

\begin{table}[htbp]
\centering
\caption{Sharp difference-regret constants and optimal quantiles.}
\label{tab:numerics}
\begin{tabular}{rcc}
\toprule
$n$ & $q_n$ & $-n\log x_n^*$\\
\midrule
2 & 0.062500000000 & 1.386294361120 \\
3 & 0.089726074648 & 1.456222612442 \\
10 & 0.139510137746 & 1.608414704579 \\
25 & 0.157237068977 & 1.675058798487 \\
50 & 0.164399487922 & 1.705415999682 \\
100 & 0.168501361559 & 1.724244242335 \\
200 & 0.170811703333 & 1.735556227083 \\
500 & 0.172371267380 & 1.743648470057 \\
1000 & 0.172953035237 & 1.746820023676 \\
2000 & 0.173269573073 & 1.748604193956 \\
3000 & 0.173382172602 & 1.749253962769 \\
10000 & 0.173551874778 & 1.750257603059 \\
100000 & 0.173627047643 & 1.750720554040 \\
\bottomrule
\end{tabular}
\end{table}

\section{Discussion}

The sharp additive comparison with the prophet is determined by two
maximizing points of the difference kernel. Their balance identifies the
optimal quantile, and their weights identify the unique least-favorable
law. This is why the continuum of distributions in the original problem
can be replaced by an explicit distribution on three observation values.
The universal rule requires only the horizon and the quantile of the
known distribution, rather than the numerical solution of a matrix game.

The behavior for large horizons is less immediate than the limiting
constant might suggest. The upper adversarial point remains within
order $1/n$ of $1$, whereas the lower point is separated from $1$
by order $(\log n)/n$. The different scales explain the logarithmic
term in the convergence rate. In particular,
Theorem~\ref{thm:asympt} implies
\[
 \frac{n(c_*-v_n)}{\log n}\longrightarrow c_*(1-\rho_*)>0.
\]
Proposition~\ref{prop:monotonicity} shows that this gap is positive at
every finite horizon, not only in the asymptotic regime. The explicit
quantile in \eqref{eq:simplex} is thus a convenient alternative when the
limiting guarantee is sufficient.

These conclusions concern the additive comparison for one fixed threshold,
with independent tie-breaking at equality and compulsory selection of the
last observation. They do not assert the same constants for adaptive
thresholds, ratio objectives, or different terminal-payoff conventions.

\section*{Acknowledgments}
During preparation of this manuscript, the author used ChatGPT (OpenAI,
GPT-5.6 Sol; September 2026; \url{https://chatgpt.com}) to assist with
exploratory derivations, language editing, literature searches, and
algebraic and numerical cross-checks.

\appendix
\section{A derivative verification of the saddle inequality}

The convexity proof of Theorem~\ref{thm:exact} establishes the global
minimum of $G_n$ defined in \eqref{eq:Gdef}. For completeness, we record
its derivative in the original quantile coordinate. This also identifies
the interval of minimizing quantiles in the exceptional two-observation
case. On the two outer intervals one has
$G_n'(x)<0$ on $(0,y_{-,n}^*)$ and $G_n'(x)>0$ on $(y_{+,n}^*,1)$.
In the middle interval,
\[
 \frac{G_n'(x)}{x^{n-2}}
 =(1-\rho_n)(n-1)y_{-,n}^*
 -\rho_n(1-y_{+,n}^*)\sum_{k=1}^{n-1}k x^{k-n+1}.
\]
Every exponent in the sum is nonpositive, so the displayed expression is
nondecreasing in $x$ and vanishes at $x_n^*$. It is strictly increasing
when $n\ge3$. Thus $G_n$ is nonincreasing up to $x_n^*$ and nondecreasing afterwards;
continuity at both breakpoints completes the global argument.
When $n=2$, the middle expression vanishes identically, giving a flat
minimum on $[1/4,3/4]$ for the least-favorable law.

\end{document}